\documentclass[a4, 10pt]{amsart}
\usepackage{amssymb}
\usepackage{amstext}
\usepackage{amsmath}
\usepackage{amscd}
\usepackage{latexsym}
\usepackage{amsfonts}

\theoremstyle{plain}
\newtheorem{thm}{Theorem}[section]
\newtheorem*{thm*}{Theorem}
\newtheorem*{cor*}{Corollary}

\newtheorem*{thma}{Theorem \ref{thma}}
\newtheorem*{cora}{Corollary \ref{cora}}
\newtheorem*{thmb}{Theorem \ref{thmb}}

\newtheorem{prop}[thm]{Proposition}
\newtheorem{lem}[thm]{Lemma}

\newtheorem{cor}[thm]{Corollary}

\theoremstyle{definition}
\newtheorem{defn}[thm]{Definition}

\newtheorem{rem}[thm]{Remark}

\newtheorem{ques}[thm]{Question}

\theoremstyle{remark}
\newtheorem*{pf}{{\sl Proof}}
\newtheorem*{apf}{{\sl Another proof of Theorem \ref{thma} (2) $\Rightarrow$ (1)}}

\numberwithin{equation}{thm}
\def\Hom{\mathrm{Hom}}

\def\Ext{\mathrm{Ext}}

\def\Coker{\mathrm{Coker}}
\def\Ker{\mathrm{Ker}}
\def\Im{\mathrm{Im}}
\def\tr{\mathrm{tr}}
\def\rank{\mathrm{rank}}

\def\m{\mathfrak m}

\def\J{\mathrm{J}}
\def\Gdim{\mathrm{G}\mathrm{dim}}
\def\GCdim{\mathrm{G}_C\mathrm{dim}}

\def\Kdim{\mathrm{dim}}
\def\depth{\mathrm{depth}}

\def\pd{\mathrm{pd}}

\def\grade{\mathrm{grade}}

\begin{document}

\title{Syzygy modules with semidualizing or G-projective summands}
\author{Ryo Takahashi}
\address{Department of Mathematics, Faculty of Science, Okayama University, 1-1, Naka 3-chome, Tsushima, Okayama 700-8530, Japan}
\email{takahasi@math.okayama-u.ac.jp}
\thanks{
\endgraf
}
\maketitle
\begin{abstract}
Let $R$ be a commutative Noetherian local ring with residue class field $k$.
In this paper, we mainly investigate direct summands of the syzygy modules of $k$.
We prove that $R$ is regular if and only if some syzygy module of $k$ has a semidualizing summand.
After that, we consider whether $R$ is Gorenstein if and only if some syzygy module of $k$ has a G-projective summand.
\end{abstract}
\section{Introduction}

Throughout the present paper, unless otherwise specified, all rings will be assumed to be commutative Noetherian local rings, and all modules will be assumed to be finitely generated modules.

The notion of a semidualizing module was introduced by Golod \cite{Golod} by the name of a suitable module.
This module is a PG-module, which was defined by Foxby \cite{Foxby} as a common generalization of a projective module and a Gorenstein module.
A free module of rank one over any local ring and the canonical module over a Cohen-Macaulay local ring are semidualizing modules.
This notion has been extended to a complex of modules, which is called a semidualizing complex, by Christensen \cite{Christensen2}.

On the other hand, the notion of a finitely generated G-projective module (in other words, a module of G-dimension zero) was introduced by Auslander \cite{Auslander} and was deeply studied by him and Bridger \cite{AB}.
All free modules over any local ring and all maximal Cohen-Macaulay modules over a Gorenstein local ring are G-projective modules.
This notion has been extended to a module which is not necessarily finitely generated by Enochs and Jenda \cite{EJ}.

In this paper, $R$ always denotes a commutative Noetherian local ring with maximal ideal $\m$ and residue class field $k$.
Dutta \cite[Corollary 1.3]{Dutta} gave the following characterization of regular local rings.
We denote by $\Omega _R ^i M$ the $i$th syzygy module of an $R$-module $M$.

\begin{thm}[Dutta]\label{duttathm}
The following are equivalent:
\begin{enumerate}
\item[$(1)$]
$R$ is regular;
\item[$(2)$]
$\Omega _R ^n k$ has a nonzero free summand for some $n\geq 0$.
\end{enumerate}
\end{thm}

This theorem suggests a natural question asking whether a certain class of local rings is characterized by the existence of a semidualizing module or a G-projective module which is isomorphic to a direct summand of some syzygy module of the residue class field.

The main purpose of this paper is to answer this question.
The first goal of this paper is the following theorem, which is a generalization of Dutta's theorem.

\begin{thma}
The following are equivalent:
\begin{enumerate}
\item[$(1)$]
$R$ is regular;
\item[$(2)$]
$\Omega _R ^n k$ has a semidualizing summand for some $n\geq 0$.
\end{enumerate}
\end{thma}

This theorem yields a criterion for a Cohen-Macaulay local ring to be regular using the canonical module:

\begin{cora}
Let $R$ be a Cohen-Macaulay local ring with canonical module $\omega$.
Then the following are equivalent:
\begin{enumerate}
\item[$(1)$]
$R$ is regular;
\item[$(2)$]
$\Omega _R ^n k$ has a summand isomorphic to $\omega$ for some $n\geq 0$.
\end{enumerate}
\end{cora}

On the other hand, the author \cite{Takahashi} proved that a local ring is Gorenstein if the second syzygy module of the residue field has a G-projective summand (more generally, a summand of finite G-dimension).
Hence we naturally expect that any local ring such that some syzygy module of the residue field has a G-projective summand will be Gorenstein.
In relation to this expectation, the second goal of this paper is the following theorem.

\begin{thmb}
The following are equivalent:
\begin{enumerate}
\item[$(1)$]
$R$ is Gorenstein;
\item[$(2)$]
$\Omega _R ^n k$ has a nonzero G-projective summand for some $0\leq n\leq\depth\,R+2$.
\end{enumerate}
\end{thmb}

The organization of this paper is as follows; this paper consists of six sections.

In Section 2, we will recall some notions which are necessary for our proofs of the main results of this paper.
The definitions of a semidualizing module and a G-projective module will be recalled in this section.

In Section 3, we will investigate lower syzygy modules.
We will prove using elementary techniques that the $i$th syzygy of the residue field $k$ of $R$ is indecomposable for any integer $i<\depth\,R$.
This result will be used in all the subsequent sections.

In Section 4, we will consider the condition that some syzygy of the residue field has a semidualizing summand, and prove Theorem \ref{thma} and Corollary \ref{cora}, which are stated above.

In Section 5, we will give another approach to Theorem \ref{thma}.
Taking into account the relationship between a semidualizing module and a nonzerodivisor, we will make another proof of the theorem.

In Section 6, we will consider whether $R$ is Gorenstein provided that some syzygy of the residue field $k$ of $R$ has a G-projective summand.
Taking advantage of a certain short exact sequence obtained from results of Auslander and Bridger, we will prove Theorem \ref{thmb}, which is stated above.

\section{Basic results on G$_C$-dimension}

In this section, we recall the definition of a semidualizing module and several notions related to it.
We also state basic properties of them which we will use in the later sections.

\begin{defn}
An $R$-module $C$ is said to be {\it semidualizing} (or {\it suitable}) if the following hold:
\begin{enumerate}
\item[$(1)$]
The natural homomorphism $R\to\Hom _R (C, C)$ is an isomorphism,
\item[$(2)$]
$\Ext _R ^i (C, C)=0$ for any $i>0$.
\end{enumerate}
\end{defn}

Note that $R$ is itself a semidualizing $R$-module.
If $R$ is Cohen-Macaulay, then the canonical module of $R$ is a semidualizing $R$-module.

\begin{prop}\label{sd}
Let $C$ be a semidualizing $R$-module.
Then the following hold.
\begin{enumerate}
\item[$(1)$]
$C$ is indecomposable.
\item[$(2)$]
$\depth _R\, C= \depth\, R$.
\item[$(3)$]
Let $x\in\m$ be an $R$-regular element.
Then $C/xC$ is a semidualizing $R/(x)$-module.
\item[$(4)$]
Let $R\to S$ be a flat homomorphism of local rings.
Then $C\otimes _R S$ is a semidualizing $S$-module.
\end{enumerate}
\end{prop}

\begin{pf}
The statements (2) and (3) are shown in \cite{Golod}.
The statement (4) is easy to see by definition.
Concerning the statement (1), suppose that there exists a direct sum decomposition $C=C'\oplus C''$ with $C', C''\not=0$.
Then we have isomorphisms $R\cong\Hom (C, C)\cong\Hom (C', C')\oplus\Hom (C', C'')\oplus\Hom (C'', C')\oplus\Hom (C'', C'')$.
The modules $\Hom (C', C')$ and $\Hom (C'', C'')$ are nonzero because they contain the identity maps of $C'$ and $C''$, respectively.
However since $R$ is a local ring, it is indecomposable as an $R$-module, and we get a contradiction.
Thus $C$ is an indecomposable $R$-module.
\qed
\end{pf}

To prove the main results of this paper, we will take advantage of the notion of G$_C$-dimension.
It is a homological invariant for modules associated to a fixed semidualizing module $C$.
As well as that of a semidualizing module, the notions of a G$_C$-projective module and G$_C$-dimension were originally introduced by Golod \cite{Golod}.

\begin{defn}
Let $C$ be a semidualizing $R$-module, and denote by $(-)^{\dag}$ the $C$-dual functor $\Hom _R (-, C)$.
\begin{enumerate}
\item[$(1)$]
We say that an $R$-module $X$ is {\it G$_C$-projective} (or {\it $C$-reflexive}) if the following hold:
\begin{enumerate}
\item[{\rm (i)}]
The natural homomorphism $X\to X^{\dag\dag}$ is an isomorphism,
\item[{\rm (ii)}]
$\Ext _R ^i (X, C)=0$ for any $i>0$,
\item[{\rm (iii)}]
$\Ext _R ^i (X^{\dag}, C)=0$ for any $i>0$.
\end{enumerate}
\item[$(2)$]
Let $M$ be an $R$-module.
If there exists an exact sequence
$$
0 \to X_n \to X_{n-1} \to \cdots \to X_1 \to X_0 \to M \to 0
$$
of $R$-modules such that each $X_i$ is G$_C$-projective, then we say that $M$ has {\it G$_C$-dimension at most $n$}, and write $\GCdim _R\, M\leq n$.
If such an integer $n$ does not exist, then we say that $M$ has {\it infinite G$_C$-dimension}, and write $\GCdim _R\, M=\infty$.
\end{enumerate}
\end{defn}

If an $R$-module $M$ has G$_C$-dimension at most $n$ but does not have G$_C$-dimension at most $n-1$, then we say that $M$ has {\it G$_C$-dimension $n$}, and write $\GCdim _R\, M=n$.
Note that being G$_C$-dimension zero is equivalent to being G$_C$-projective.
Also note that $C$ is itself G$_C$-projective, and so is any free $R$-module.

G$_C$-dimension has similar properties to projective dimension:

\begin{prop}\label{gd}
Let $C$ be a semidualizing $R$-module.
Denote by $(-)^{\dag}$ the $C$-dual functor $\Hom _R (-, C)$.
\begin{enumerate}
\item[$(1)$]\ Let $M$ be a nonzero $R$-module with $\GCdim _R\, M<\infty $.
Then
\begin{align*}
\GCdim _R\, M & = \depth\, R-\depth _R\, M \\
& = \sup\{\,i\,\vert\,\Ext _R ^i (M, C)\not =0\,\}.
\end{align*}
\item[$(2)$]\ Let $0 \to L \to M \to N \to 0$ be an exact sequence of $R$-modules.
If two of $L, M, N$ have finite G$_C$-dimension, then so does the third.
\item[$(3)$]\ Let $M$ be an $R$-module.
Then
$$
\GCdim _R (\Omega _R ^n M)=\sup\{\,\GCdim _R\, M-n,\,0\,\}
$$
for any $n\geq 0$.
\item[$(4)$]\ Let $M, N$ be $R$-modules.
Then
$$
\GCdim _R (M\oplus N)=\sup\{\,\GCdim _R\, M,\,\GCdim _R\, N\,\}.
$$
\end{enumerate}
\end{prop}

\begin{pf}
The statements (1) and (2) are stated in \cite{Golod}, and the other statements are easily observed from them.
\qed
\end{pf}

A {\it G-projective} $R$-module is defined as a G$_R$-projective $R$-module.
The {\it G-dimension} of an $R$-module $M$ is defined as the G$_R$-dimension of $M$, and is denoted by $\Gdim _R\, M$.

We denote by $(-)^{\ast}$ the $R$-dual functor $\Hom _R (-, R)$.
Let
$$
F_1 \overset{\delta}{\to} F_0 \to M \to 0
$$
be a free presentation of an $R$-module $M$.
We define by $\tr _R M$ the cokernel of the $R$-dual homomorphism ${\delta}^{\ast} : F_0 ^{\ast} \to F_1 ^{\ast}$.
It is called the {\it Auslander transpose} or the {\it Auslander dual} of $M$.
The module $\tr _R M$ depends on the choice of a free presentation of $M$, but is uniquely determined up to free summand.
For more details, refer to \cite{AB} and \cite{Masek}.

In the last section of this paper we will need the following result concerning G-projective modules and G-dimension.

\begin{prop}\label{orth}
\begin{enumerate}
\item[$(1)$]
If an $R$-module $M$ is G-projective, then so are $M^{\ast}$, $\Omega _R M$ and $\tr _R M$.
\item[$(2)$]
Let $X$ be a G-projective $R$-module, and let $T$ be an $R$-module of finite projective dimension.
Then $\Gdim _R\,\Hom _R (X, T)<\infty$ and $\Ext _R ^i (X, T)=0$ for any $i>0$.
\item[$(3)$]
The following are equivalent:
\begin{enumerate}
\item[{\rm (i)}]
$R$ is Gorenstein;
\item[{\rm (ii)}]
$\Gdim _R\, M<\infty $ for any $R$-module $M$;
\item[{\rm (iii)}]
$\Gdim _R\, k<\infty $.
\end{enumerate}
\end{enumerate}
\end{prop}

\begin{pf}
We can find proofs of the assertions (2) and (3) in \cite[Proposition (4.12)]{AB} and \cite[(1.4.9)]{Christensen} respectively.
For the assertion (1), we easily observe from definition that $M^{\ast}$ is G-projective.
It is seen by Proposition \ref{gd}(3) that $\Omega M$ is G-projective.
That $\tr M$ is G-projective follows from \cite[Proposition (3.8)]{AB}.
\qed
\end{pf}

\section{Indecomposability of syzygy modules}

We study in this section the indecomposability of syzygy modules in elementary ways.
First of all, we give an easy lemma.

\begin{lem}\label{five}
Let
$$
0 \to L \overset{f}{\to} M \overset{g}{\to} N \to 0
$$
be an exact sequence of $R$-modules.
Suppose that $L$ has a direct sum decomposition $L=X\oplus Y$, and write $f = (\sigma, \tau)$ along this decomposition.
Then we have the following exact sequences of $R$-modules:
$$
\begin{cases}
0 \longrightarrow X \overset{\sigma}{\longrightarrow} M \overset{\alpha}{\longrightarrow} A \longrightarrow 0, & \quad\cdots (1) \\
0 \longrightarrow Y \overset{\tau}{\longrightarrow} M \overset{\beta}{\longrightarrow} B \longrightarrow 0, & \quad\cdots (2) \\
0 \longrightarrow X \overset{\beta\sigma}{\longrightarrow} B \overset{\zeta}{\longrightarrow} N \longrightarrow 0, & \quad\cdots (3) \\
0 \longrightarrow Y \overset{\alpha\tau}{\longrightarrow} A \overset{\eta}{\longrightarrow} N \longrightarrow 0, & \quad\cdots (4) \\
0 \longrightarrow M \overset{\binom{\alpha}{\beta}}{\longrightarrow} A\oplus B \overset{\gamma}{\longrightarrow} N \longrightarrow 0. & \quad\cdots (5)
\end{cases}
$$
\end{lem}

\begin{pf}
It is easy to get the exact sequences (1)--(4).
Setting $A = \Coker\,\sigma$ and $B = \Coker\,\tau$, we get (1) and (2).
Defining the map $\zeta$ by $\zeta (\beta (x))= g(x)$ for $x\in M$ gives (3).
Similarly we obtain (4).

As for the exact sequence (5), we define the map $\gamma$ by $\gamma (\binom{\alpha (x)}{\beta (y)})=g(x-y)$ for $x, y\in M$.
If $\gamma (\binom{\alpha (x)}{\beta (y)})=g(x-y)=0$, then we have $x-y=f(\binom{s}{t})=\sigma (s)+\tau (t)$ for some $\binom{s}{t}\in L$.
Putting $z=x-\sigma (s)=y+\tau (t)$, we have $\binom{\alpha (x)}{\beta (y)}=\binom{\alpha}{\beta}(z)$.
Hence $\Im\binom{\alpha}{\beta}=\Ker\,\gamma$, and thus (5) is obtained.
\qed
\end{pf}

We see from the proof that the above result holds for an arbitrary ring $R$ and arbitrary $R$-modules $L, M, N$.

For an $R$-module $M$, we denote by $\nu _R (M)$ the minimal number of generators of $M$, i.e., $\nu _R (M)= \Kdim _k (M\otimes _R k)$.

\begin{prop}\label{grade}
Suppose that $R$ is Henselian.
Let $M$ be an indecomposable $R$-module.
Then $\Omega _R ^i M$ is indecomposable for any $i$ with $0\leq i<\grade _R\, M$.
\end{prop}

\begin{pf}
Set $t=\grade\, M$.
We show the proposition by induction on $i$.
It trivially holds for $i=0$, and so let $i\geq 1$.
Take a minimal free resolution $F_{\bullet}$ of $M$.
Suppose that $\Omega ^i M$ has a direct sum decomposition $\Omega ^i M = X \oplus Y$.
Then we have an exact sequence
$$
0 \to X\oplus Y \to F_{i-1} \to \Omega ^{i-1} M \to 0.
$$
According to Lemma \ref{five}, there exist exact sequences
\begin{align}
\label{grade1}
& 0 \to X \to F_{i-1} \to A \to 0,\\
\label{grade2}
& 0 \to Y \to F_{i-1} \to B \to 0,\\
\label{grade3}
& 0 \to X \to B \to \Omega ^{i-1} M \to 0,\\
\label{grade4}
& 0 \to Y \to A \to \Omega ^{i-1} M \to 0,\\
\label{grade5}
& 0 \to F_{i-1} \to A\oplus B \to \Omega ^{i-1} M \to 0.
\end{align}
The last exact sequence \eqref{grade5} can be regarded as an element of $\Ext ^1 (\Omega ^{i-1} M, F_{i-1})$, and we have $\Ext ^1 (\Omega ^{i-1} M, F_{i-1})\cong\Ext ^i (M, F_{i-1})=0$ because $i<t$.
Hence the exact sequence \eqref{grade5} splits; we have an isomorphism
$$
A \oplus B \cong F_{i-1} \oplus \Omega ^{i-1} M.
$$
Since $\Omega ^{i-1} M$ is indecomposable by the induction hypothesis, it is isomorphic to a direct summand of $A$ or $B$ by the Krull-Schmidt theorem.
We can assume without loss of generality that $\Omega ^{i-1} M$ is isomorphic to a direct summand of $A$.
Then $B$ is isomorphic to a direct summand of the free module $F_{i-1}$, hence $B$ is a free module.
Therefore the exact sequence \eqref{grade2} splits, and so $Y$ is also a free module with $\rank\,Y=\rank\,F_{i-1}-\rank\,B$.
Since there is a surjective homomorphism from $B$ to $\Omega ^{i-1} M$ by \eqref{grade3}, we have
$$
\rank\,B=\nu _R (B)\geq\nu _R (\Omega ^{i-1}M)=\rank\,F_{i-1}.
$$
Thus we have $\rank\,F_{i-1}=\rank\,B$, and $Y=0$.
This means that the $R$-module $\Omega ^i M$ is indecomposable.
\qed
\end{pf}

As a corollary of the above result, we obtain the following result.
We should note that $R$ is not assumed to be Henselian.

\begin{cor}\label{below}
$\Omega _R ^i k$ is indecomposable for $0\leq i<\depth\, R$.
\end{cor}

\begin{pf}
Fix an integer $i$ with $0\leq i<\depth\, R$.
Let $\Omega _R ^i k = X\oplus Y$ be a direct sum decomposition.
Taking the $\m$-adic completions, we obtain $\Omega _{\widehat{R}} ^i k \cong \widehat{X}\oplus\widehat{Y}$.
Since $\widehat{R}$ is Henselian and $k$ is an indecomposable $\widehat{R}$-module, $\Omega _{\widehat{R}} ^i k$ is also an indecomposable $\widehat{R}$-module by Proposition \ref{grade}, and we have $\widehat{X}=0$ or $\widehat{Y}=0$.
Hence $X=0$ or $Y=0$, which says that $\Omega _R ^i k$ is indecomposable as an $R$-module.
\qed
\end{pf}

\section{Syzygy modules with semidualizing summands}

In the previous section we demonstrated that the $i$th syzygy module of the residue field $k$ is indecomposable for $i<\depth\,R$.
We start this section by showing that any $i$th syzygy module has no nonfree semidualizing summand for $i>\depth\,R$.

\begin{prop}\label{above}
Let $M$ be an $R$-module and $C$ a semidualizing $R$-module.
If $C$ is isomorphic to a direct summand of $\Omega _R ^i M$ for some $i>\depth\,R$, then $C\cong R$.
\end{prop}

\begin{pf}
Fix an integer $i>t:=\depth\,R$, and take a minimal free resolution $F_{\bullet}$ of $M$.
There is an isomorphism $\Omega ^i M\cong C\oplus D$ for some $R$-module $D$.
We have an exact sequence
$$
0 \to C\oplus D \to F_{i-1} \to \Omega ^{i-1}M \to 0.
$$
It follows from Lemma \ref{five} that we have exact sequences
\begin{align}
\label{above1}
& 0 \to C \to F_{i-1} \to A \to 0,\\
\label{above2}
& 0 \to D \to F_{i-1} \to B \to 0,\\
\label{above3}
& 0 \to C \to B \to \Omega ^{i-1}M \to 0,\\
\label{above4}
& 0 \to D \to A \to \Omega ^{i-1}M \to 0.
\end{align}
From the exact sequence $\eqref{above1}$ and Proposition \ref{gd}(2) and \ref{gd}(1), we see that $\GCdim\,A=\sup\{\,i\,\vert\,\Ext ^i (A, C)\not=0\,\}\leq 1$.

Suppose that $\Ext ^1 (A, C)\not=0$.
Then $\GCdim\,A=1$, and $\depth\,A=t-1$ by Proposition \ref{gd}(1).
We have
$$
\depth\,\Omega ^{i-1}M\geq\inf\{ i-1, \depth\,R\} =t
$$
by \cite[Exercise 1.3.7]{BH},
$$
\depth\,B\geq\inf\{\depth\,C, \depth\,\Omega ^{i-1}M\} =t
$$
by \eqref{above3} and Proposition \ref{sd}(2), and
$$
\depth\,D\geq\inf\{\depth\,F_{i-1}, \depth\,B+1\} =t
$$
by \eqref{above2}.
Therefore by \eqref{above4} we get
$$
t-1=\depth\,A\geq\inf\{\depth\,D, \depth\,\Omega ^{i-1}M\}\geq t,
$$
which is a contradiction.

Thus we must have $\Ext ^1 (A, C)=0$.
Hence the exact sequence $\eqref{above1}$ splits, therefore $C$ is free, and thus $C\cong R$ by Proposition \ref{sd}(1).
\qed
\end{pf}

For $i=\depth\,R$, the $i$th syzygy module of the residue field $k$ has no nonfree semidualizing summand:

\begin{prop}\label{just}
If $C$ is a semidualizing $R$-module which is isomorphic to a direct summand of $\Omega _R ^t k$ where $t=\depth\,R$, then $C\cong R$.
\end{prop}

\begin{pf}
When $t=0$, we have $C\cong k$ since $k$ is an indecomposable $R$-module.
Hence $R\cong\Hom _R (C, C)=\Hom _R (k, k)\cong k\cong C$.
Let $t\geq 1$ and take a minimal free resolution $F_{\bullet}$ of the $R$-module $k$.
We have an exact sequence $0 \to \Omega ^t k \to F_{t-1} \to \Omega ^{t-1} k \to 0$.
There is an isomorphism $\Omega ^t k\cong C\oplus D$ for some $R$-module $D$.
Then Lemma \ref{five} provides exact sequences
\begin{align}
\label{just1}
& 0 \longrightarrow C \overset{\sigma}{\longrightarrow} F_{t-1} \overset{\alpha}{\longrightarrow} A \longrightarrow 0,\\
\label{just2}
& 0 \longrightarrow D \overset{\tau}{\longrightarrow} F_{t-1} \overset{\beta}{\longrightarrow} B \longrightarrow 0,\\
\label{just3}
& 0 \longrightarrow C \overset{\beta\sigma}{\longrightarrow} B \longrightarrow \Omega ^{t-1} k \longrightarrow 0,\\
\label{just4}
& 0 \longrightarrow D \overset{\alpha\tau}{\longrightarrow} A \longrightarrow \Omega ^{t-1} k \longrightarrow 0,\\
\label{just5}
& 0 \longrightarrow F_{t-1} \overset{\binom{\alpha}{\beta}}{\longrightarrow} A\oplus B \longrightarrow \Omega ^{t-1} k \longrightarrow 0.
\end{align}

Denote by $(-)^{\dagger}$ the $C$-dual functor $\Hom _R (-, C)$.
From \eqref{just3} we have an exact sequence
$$
0 \to (\Omega ^{t-1}k)^{\dagger} \to B^{\dagger} \to R \overset{\lambda}{\to} \Ext ^1 (\Omega ^{t-1}k, C) \to \Ext ^1 (B, C) \to 0.
$$
Put $I=\Ker\,\lambda$.
The map $\lambda$ sends $1\in R$ to an element of $\Ext ^1 (\Omega ^{t-1}k, C)$ corresponding to the exact sequence \eqref{just3}.
If the exact sequence \eqref{just3} splits, then the map $\beta\sigma$ is a split monomorphism, and so is the map $\sigma$.
Hence $C$ is free, and therefore $C\cong R$ by Proposition \ref{sd}(1).
Thus we may assume that \eqref{just3} does not split.
Then the map $\lambda$ in the above exact sequence is nonzero, equivalently $I\not=R$.
Noting that $\Ext ^1 (\Omega ^{t-1}k, C)$ is a $k$-vector space, we have $I=\m$, and
\begin{equation}\label{vect}
\Kdim _k\, \Ext ^1 (B, C)=\Kdim _k\, \Ext ^1 (\Omega ^{t-1}k, C)-1.
\end{equation}

From \eqref{just1} we have an exact sequence
$$
0 \to A^{\dagger} \to F_{t-1}^{\dagger} \overset{\rho}{\to} R \to \Ext ^1 (A, C) \to 0.
$$
Put $J=\Im\,\rho$.
If $\Ext ^1 (A, C)=0$, then the exact sequence \eqref{just1} splits, $C$ is free, and $C\cong R$ by Proposition \ref{sd}(1).
Hence we may assume that $\Ext ^1 (A, C)\not=0$, equivalently $J\not=R$.
From \eqref{just5} we have an exact sequence
\begin{align*}
0 & \to (\Omega ^{t-1}k)^{\dagger} \to A^{\dagger}\oplus B^{\dagger} \to F_{t-1}^{\dagger} \\
& \to \Ext ^1 (\Omega ^{t-1}k, C) \overset{\mu}{\to} \Ext ^1 (A, C)\oplus\Ext ^1 (B, C) \to 0.
\end{align*}
Hence $\Ext ^1 (A, C)$ is a $k$-vector space since this is true of $\Ext ^1 (\Omega ^{t-1}k, C)$.
Therefore we see that $J=\m$ and $\Kdim _k\,\Ext ^1 (A, C)=1$.
It follows from \eqref{vect} that $\mu$ must be an isomorphism.
Thus we obtain two exact sequences
\begin{gather}
\label{adag}
0 \to A^{\dagger} \to F_{t-1}^{\dagger} \to \m \to 0,\\
\label{abdag}
0 \to (\Omega ^{t-1}k)^{\dagger} \to A^{\dagger}\oplus B^{\dagger} \to F_{t-1}^{\dagger} \to 0.
\end{gather}

Note that $\depth\,C=t\geq 1$ by Proposition \ref{sd}(2).
Applying $(-)^{\dagger}$ to the complex $F_{\bullet}$ gives an exact sequence
$$
0 \to F_0^{\dagger} \to \dots \to F_{t-2}^{\dagger} \to (\Omega ^{t-1}k)^{\dagger} \to 0.
$$
Noting that each $F_i^{\dagger}$ is a direct sum of copies of $C$ and applying the functor $\Hom (C, -)$ to this exact sequence, we easily see that $\Ext ^i (C, (\Omega ^{t-1}k)^{\dagger})=0$ for any $i>0$.
Hence applying $\Hom (C, -)$ to \eqref{abdag}, we have $\Ext ^i (C, A^{\dagger})=0$ for any $i>0$, and applying $\Hom (C, -)$ to \eqref{adag}, we have $\Ext ^i (C, \m)=0$ for any $i>0$.
In particular, $\Ext ^1 (C, \m)=0$, which shows that any homomorphism from $C$ to $k$ factors through the natural surjection $\pi : R \to k$.
Choose a minimal generator $x$ of $C$ and a homomorphism $\varepsilon : C \to k$ such that $\varepsilon (x)=\pi (1)$.
There exists a homomorphism $\xi : C \to R$ such that $\varepsilon = \pi\xi$.
Then it is easy to see that $\xi (x)$ is a unit in $R$, which implies that $\xi$ is surjective.
Hence $\xi$ is a split epimorphism, and thus $C\cong R$ by Proposition \ref{sd}(1).
\qed
\end{pf}

Applying Dutta's theorem to our results obtained above, we give a characterization of regular local rings, which is one of the main results of this paper.

\begin{thm}\label{thma}
The following are equivalent:
\begin{enumerate}
\item[$(1)$]
$R$ is regular;
\item[$(2)$]
$\Omega _R ^n k$ has a semidualizing summand for some $n\geq 0$.
\end{enumerate}
\end{thm}

\begin{pf}
(1) $\Rightarrow$ (2): Set $d=\Kdim\,R$.
Then the $R$-module $\Omega ^d k$ is free of rank one.
Hence it is itself semidualizing.

(2) $\Rightarrow$ (1): Assume $n<\depth\,R$.
Then it follows from Corollary \ref{below} that $C$ is isomorphic to $\Omega ^n k$.
We easily see that $\depth\,\Omega ^n k=n$, while we have $\depth\,C=\depth\,R>n$ by Proposition \ref{sd}(2), which is a contradiction.
Therefore $n\geq\depth\,R$.
Proposition \ref{above} and \ref{just} say that $C$ is isomorphic to $R$, and Theorem \ref{duttathm} shows that $R$ is regular.
\qed
\end{pf}

As an application of the above theorem, we obtain a necessary and sufficient condition for a Cohen-Macaulay local ring to be regular, using the canonical module.

\begin{cor}\label{cora}
Let $R$ be a Cohen-Macaulay local ring with canonical module $\omega$.
Then the following are equivalent:
\begin{enumerate}
\item[$(1)$]
$R$ is regular;
\item[$(2)$]
$\Omega _R ^n k$ has a summand isomorphic to $\omega$ for some $n\geq 0$.
\end{enumerate}
\end{cor}

\begin{pf}
The implication (1) $\Rightarrow$ (2) has been proved in the proof of Theorem \ref{thma}.
The converse immediately follows from Theorem \ref{thma} and the fact that the canonical module of a Cohen-Macaulay local ring is semidualizing.
\qed
\end{pf}

\section{Another approach to Theorem \ref{thma}}

In this section, we consider Theorem \ref{thma} from another point of view.
To be concrete, we direct our attention to the relationship between a semidualizing module and a regular element which is stated in Proposition \ref{sd}(3), and give a proof of the theorem which is different from the proof given in the previous section.

Let us recall the notion of a weak lifting:

\begin{defn}
Let $x\in\m$ be an $R$-regular element and set $\overline{(-)}=(-)\otimes _R R/(x)$.
An $\overline{R}$-module $M$ is said to be {\it weakly liftable} to $R$ if there exists an $R$-module $N$ such that $x$ is $N$-regular and $M$ is a direct summand of $\overline{N}$.
\end{defn}

\begin{prop}\label{wl}
Let $x\in\m$ be $R$-regular and set $\overline{(-)}=(-)\otimes _R R/(x)$.
The following are equivalent for an $\overline{R}$-module $M$:
\begin{enumerate}
\item[$(1)$]
$M$ is weakly liftable to $R$;
\item[$(2)$]
$\overline{\Omega _R M}\cong\Omega _{\overline{R}}M\oplus M$;
\item[$(3)$]
$\overline{\Omega _R ^{i+1}M}\cong\Omega _{\overline{R}}^{i+1}M\oplus\Omega _{\overline{R}}^iM$ for any $i\geq 0$.
\end{enumerate}
\end{prop}

\begin{pf}
(1) $\Leftrightarrow$ (2): We refer to \cite[Proposition 3.2]{ADS}.

(3) $\Rightarrow$ (2): This implication is trivial.

(2) $\Rightarrow$ (3): We have an exact sequence
$$
0 \to \Omega _R ^2 M \overset{f}{\to} F \to \Omega _R M \to 0,
$$
where $F$ is a free $R$-module and $f\otimes _R k=0$.
Since $x$ is $R$-regular, it is also $\Omega _R M$-regular.
Hence we see from the snake lemma that the above exact sequence induces an exact sequence
$$
0 \to \overline{\Omega _R ^2 M} \overset{\overline{f}}{\to} \overline{F} \to \overline{\Omega _R M} \to 0,
$$
and $\overline{f}\otimes _{\overline{R}} k=0$.
Therefore we have $\overline{\Omega _R ^2 M}\cong\Omega _{\overline{R}} (\overline{\Omega _R M})\cong\Omega _{\overline{R}}(\Omega _{\overline{R}} M\oplus M)\cong\Omega _{\overline{R}}^2M\oplus\Omega _{\overline{R}}M$.
Repeating this argument, we see that the condition (3) holds.
\qed
\end{pf}

The above result yields the structure of syzygy modules of the residue class field modulo a regular element.

\begin{cor}\label{m^2}
Let $x\in\m - \m ^2$ be an $R$-regular element.
Put $\overline{(-)}=(-)\otimes _R R/(x)$.
Then
$$
\overline{\Omega _R ^{i+1} k}\cong\Omega _{\overline{R}}^{i+1}k\oplus\Omega _{\overline{R}}^ik
$$
for any $i\geq 0$.
\end{cor}

\begin{pf}
It is easily checked that the natural exact sequence
$$
0 \to k \to \m /x\m \to \m /xR \to 0
$$
splits, which implies that the $\overline{R}$-module $k$ is weakly liftable to $R$.
Thus the assertion of the corollary follows from Proposition \ref{wl}.
\qed
\end{pf}

Now, we can achieve the aim of this section.

\begin{apf}
Let $C$ be a semidualizing $R$-module which is a direct summand of $\Omega _R ^n k$.
Then by Proposition \ref{sd}(4), $\widehat{C}$ is a semidualizing $\widehat{R}$-module which is a direct summand of $\widehat{\Omega _R ^n k}\cong\Omega _{\widehat{R}} ^n k$, where $\widehat{(-)}$ denotes the $\m$-adic completion.
Replacing $R$ by $\widehat{R}$, we may assume that $R$ is complete.
We easily see $n\geq\depth\,R$ by Proposition \ref{sd}(2), Corollary \ref{below} and the fact that $\depth\,\Omega _R ^i k=i$ for $i<\depth\,R$.

Suppose $\depth\,R>0$.
Then we can choose an $R$-regular element $x\in\m - \m ^2$.
Setting $\overline{(-)}=(-)\otimes _R R/(x)$, we see from Proposition \ref{sd}(3) and Corollary \ref{m^2} that $\overline{C}$ is a semidualizing $\overline{R}$-module which is a direct summand of $\overline{\Omega _R ^n k}\cong\Omega _{\overline{R}} ^n k \oplus \Omega _{\overline{R}} ^{n-1} k$.
Since $\overline{C}$ is indecomposable by Proposition \ref{sd}(1), it follows from the Krull-Schmidt theorem that $\overline{C}$ is isomorphic to a direct summand of $\Omega _{\overline{R}} ^n k$ or $\Omega _{\overline{R}} ^{n-1} k$.
By an inductive argument, we may assume $\depth\,R=0$.

Note then that $\Omega _R ^n k$ is a submodule of a free $R$-module (even if $n=0$), hence so is $C$.
Thus we have an exact sequence of this form
$$
0 \to C \to F \to L \to 0,
$$
where $F$ is a free $R$-module.
Proposition \ref{gd}(2) implies that the module $L$ has finite G$_C$-dimension, and Proposition \ref{gd}(1) yields equalities
$$
\sup\{\,i\,\vert\,\Ext ^i (L, C)\not=0\,\} = \GCdim\,L =\depth\,R -\depth\,L=0.
$$
In particular we have $\Ext ^1 (L, C)=0$.
Hence the above exact sequence splits, and thus $C$ is free.
Theorem \ref{duttathm} shows that $R$ is regular.
\qed
\end{apf}

\section{Syzygy modules with G-projective summands}

In this section, which is the last section of this paper, we shall aim at G-projective summands of syzygy modules, and consider whether the condition that some syzygy module of the residue class field $k$ of $R$ has a G-projective summand implies the Gorensteinness of the local ring $R$.

Let $M$ be an $R$-module and $n$ an integer.
We say that $M$ is {\it $n$-torsionfree} if $\Ext _R ^i (\tr M, R)=0$ for any $1\leq i\leq n$.
Here we state a criterion for the torsionfree property.
Note by definition that the grade of the zero module is infinity; hence it is greater than all integers.

\begin{lem}\cite[Proposition (2.26)]{AB}\label{tf}
Let $M$ be an $R$-module.
The following are equivalent for an integer $n$:
\begin{enumerate}
\item[$(1)$]
$\Omega _R ^i M$ is $i$-torsionfree for $1\leq i\leq n$;
\item[$(2)$]
$\grade _R\,\Ext _R ^i (M, R)\geq i-1$ for $1\leq i\leq n$.
\end{enumerate}
\end{lem}

For a nonnegative integer $n$, we denote by $\J _n$ the composite functor $\tr\cdot\Omega ^n$, and set $\J _n ^2 = \J _n\cdot\J _n$.
The following lemma is shown in the proof of \cite[Proposition (2.21)]{AB}.

\begin{lem}\label{j}
Let $M$ be an $R$-module and let $n\geq 1$.
If $\Omega _R ^n M$ is $n$-torsionfree, then there exists an exact sequence
$$
0 \to T \to {\J}_n^2M \oplus R^s \to M \to 0
$$
of $R$-modules with $\pd _R\,T<n$.
\end{lem}

Making use of the above two lemmas, we obtain a special type of exact sequence:

\begin{prop}\label{key}
For each $i$ with $1\leq i\leq \depth\,R+1$, there exists an exact sequence
$$
0 \to T_i \to \tr\Omega ^i \tr (\Omega ^{i+1} k) \oplus R^{s_i} \to \m \to 0
$$
of $R$-modules with $\pd _R\,T_i<i$.
\end{prop}

\begin{pf}
Put $t=\depth\,R$.
The module $\Ext ^i (\m , R)$ vanishes for $1\leq i\leq t-2$, and is a $k$-vector space for $i\geq t-1$.
Hence
$$
\grade\,\Ext ^i (\m , R)
\begin{cases}
=\infty & \text{for }1\leq i\leq t-2,\\
\geq t & \text{for }i\geq t-1.
\end{cases}
$$
Therefore we have $\grade\,\Ext ^i (\m , R)\geq i-1$ for $1\leq i\leq t+1$, and Lemma \ref{tf} implies that $\Omega ^i \m$ is $i$-torsionfree for $1\leq i\leq t+1$.
Thus Lemma \ref{j} completes the proof.
\qed
\end{pf}

\begin{rem}\label{m}
Let $X$ be a nonfree indecomposable G-projective $R$-module.
Then there are isomorphisms $\Hom _R (X, \m )\cong X^{\ast}$ and $\Ext _R ^1 (X, \m )\cong\Hom _R (X, k)$.

Indeed, applying the functor $\Hom (X, -)$ to the natural exact sequence $0 \to \m \to R \to k \to 0$, we get another exact sequence
$$
0 \to \Hom (X, \m) \overset{f}{\to} X^{\ast} \overset{g}{\to} \Hom (X, k) \overset{h}{\to} \Ext ^1 (X, \m) \to \Ext ^1 (X, R)=0.
$$
Since $X$ is nonfree and indecomposable, any homomorphism $X$ to $R$ factors through $\m$.
Hence $f$ is an isomorphism, $g$ is the zero map, and $h$ is an isomorphism.
\end{rem}

We have reached the stage to prove one of the main results of this paper.
The exact sequence given in Proposition \ref{key} plays an essential role in the proof.

\begin{thm}\label{thmb}
The following are equivalent:
\begin{enumerate}
\item[$(1)$]
$R$ is Gorenstein;
\item[$(2)$]
$\Omega _R ^n k$ has a nonzero G-projective summand for some $0\leq n\leq\depth\,R+2$.
\end{enumerate}
\end{thm}

\begin{pf}
(1) $\Rightarrow$ (2): Proposition \ref{orth}(3), Proposition \ref{gd}(1) and \ref{gd}(3) say that $\Omega ^n k$ is itself G-projective, where $n=\depth\,R\,(\leq\depth\,R+2)$.

(2) $\Rightarrow$ (1): Put $t=\depth\,R$.
Let $X$ be a nonzero G-projective summand of $\Omega ^n k$.
Then $\Omega ^{t+2-n} X$ is a G-projective summand of $\Omega ^{t+2} k$ by Proposition \ref{orth}(1).
If $\Omega ^{t+2-n} X=0$, then $\pd\,X<\infty$, and we have $\pd\,X=\depth\,R-\depth\,X=\Gdim\,X=0$ by Proposition \ref{gd}(1), hence $X$ is free, and $R$ is regular by Theorem \ref{duttathm}.
Therefore we have only to consider when $\Omega ^{t+2-n} X\not=0$.
Replacing $X$ by $\Omega ^{t+2-n} X$, we may assume $n=t+2$.

Set $E=\tr\Omega ^{t+1} \tr (\Omega ^{t+2} k)$ and $Y=\tr\Omega ^{t+1} \tr X$.
By Proposition \ref{key} we have an exact sequence
\begin{equation}\label{keyses}
0 \to T \to E\oplus R^s \to \m \to 0
\end{equation}
with $\pd\,T<\infty$.

Fix a nonfree indecomposable G-projective $R$-module $U$.
Apply the functor $\Hom (U, -)$ to \eqref{keyses}.
From Proposition \ref{orth}(2) and Remark \ref{m}, we get an exact sequence
$$
0 \to \Hom (U, T) \to \Hom (U, E)\oplus (U^{\ast})^s \to U^{\ast} \to 0,
$$
and isomorphisms
\begin{equation}\label{isom}
\Ext ^i (U, E)\cong\Ext ^i(U, \m )
\end{equation}
for $i>0$, and $\Gdim\,\Hom (U, T)<\infty$.
It is seen by Propositions \ref{orth}(1), \ref{gd}(2) and \ref{gd}(4) that $\Gdim\,\Hom (U, E)<\infty$.
Since $X$ is a direct summand of $\Omega ^{t+2} k$, the $R$-module $Y$ is a direct summand of $E$, and $\Hom (U, Y)$ is a direct summand of $\Hom (U, E)$.
Hence, by Proposition \ref{gd}(4), we have
\begin{equation}\label{why}
\Gdim\,\Hom (U, Y)<\infty\text{ for any G-projective }R\text{-module }U.
\end{equation}
Note from Proposition \ref{orth}(1) that $Y$ is itself G-projective.

Fix any G-projective $R$-module $Z$.
There is an exact sequence $0 \to \Omega Z \to F \to Z \to 0$ of $R$-modules, where $F$ is free.
Dualizing this sequence by $Y$, we obtain an exact sequence
$$
0 \to \Hom (Z, Y) \to \Hom (F, Y) \to \Hom (\Omega Z, Y) \to \Ext ^1 (Z, Y) \to 0.
$$
Since the $R$-modules $Z$, $F$, $\Omega Z$ are G-projective by Proposition \ref{orth}(1), we see from \eqref{why} that the $R$-modules $\Hom (Z, Y)$, $\Hom (F, Y)$, $\Hom (\Omega Z, Y)$ are of finite G-dimension.
Therefore, by Proposition \ref{gd}(2), we easily see that
\begin{equation}\label{gext}
\Gdim\,\Ext ^1 (Z, Y)<\infty\text{ for any G-projective }R\text{-module }Z.
\end{equation}
On the other hand, it is observed from \eqref{isom} that there is an isomorphism $\Ext ^1 (Z, E)\cong\Ext ^1 (Z, \m )$, and from Remark \ref{m} that $\Ext ^1 (Z, \m )$ is a $k$-vector space.
Hence $\Ext ^1 (Z, E)$ is also a $k$-vector space.
Since $\Ext ^1 (Z, Y)$ is a direct summand of $\Ext ^1 (Z, E)$,
\begin{equation}\label{extk}
\Ext ^1 (Z, Y)\text{ is a }k\text{-vector space for any G-projective }R\text{-module }Z.
\end{equation}

Now, assume that $R$ is not Gorenstein.
It then follows from \eqref{gext}, \eqref{extk}, Proposition \ref{gd}(4) and \ref{orth}(3) that we must have
$$
\Ext ^1 (Z, Y)=0\text{ for any G-projective }R\text{-module }Z.
$$
There is an exact sequence $0 \to \Omega (Y^{\ast}) \to P \to Y^{\ast} \to 0$ of $R$-modules where $P$ is free, and dualizing this, we get another exact sequence
$$
0 \to Y \to P^{\ast} \to (\Omega (Y^{\ast}))^{\ast} \to 0.
$$
Since $(\Omega (Y^{\ast}))^{\ast}$ is G-projective by Proposition \ref{orth}(1), we have $\Ext ^1 ((\Omega (Y^{\ast}))^{\ast}, Y)=0$, which implies that the above exact sequence splits.
Hence $Y$ is free, and so is $\tr Y$.
It follows from the definition of $Y$ that $\Omega ^{t+1} \tr X$ is also free.
Therefore we have $\pd (\tr X)<\infty$, and $\pd (\tr X)=\depth\,R-\depth (\tr X)=\Gdim (\tr X)=0$ by Proposition \ref{gd}(1) and \ref{orth}(1), that is to say, $\tr X$ is free, and so is $X$.
Thus Theorem \ref{duttathm} shows that $R$ is regular, contrary to our assumption that $R$ is not Gorenstein.
This contradiction proves the theorem.
\qed
\end{pf}

We close this paper by presenting a natural question.

\begin{ques}
Suppose that $\Omega _R ^n k$ has a nonzero G-projective summand for some integer $n>\depth\,R+2$.
Then is $R$ Gorenstein?
\end{ques}

{\sc Acknowledgments.}
The author would like to thank the referee for his/her careful reading and useful suggestions.


\end{document}